# GEOMETRY OF SOME TAXICAB CURVES

Maja Petrović [1]
Branko Malešević [2]
Bojan Banjac [3]
Ratko Obradović [4]

**Abstract**

*In this paper we present geometry of some curves in Taxicab metric. All curves of second order and trifocal ellipse in this metric are presented. Area and perimeter of some curves are also defined.*

***Key words:*** *Taxicab metric, Conics, Trifocal ellipse*

## 1. INTRODUCTION

In this paper, taxicab and standard Euclidean metrics for a visual representation of some planar curves are considered. Besides the term taxicab, also used are Manhattan, rectangular metric and city block distance [4], [5], [7]. This metric is a special case of the Minkowski metrics of order $k$ $(k \geq 1)$, which is for distance between two points $A(x_A, y_A)$ and $B(x_B, y_B)$ determined by:

$$d_k(A,B) = \left(|x_A - x_B|^k + |y_A - y_B|^k\right)^{\frac{1}{k}} \quad (1)$$

---

[1] MSc Maja Petrović, assistant at Faculty of Transport and Traffic Engineering, University of Belgrade, Serbia, e-mail: majapet@sf.bg.ac.rs

[2] PhD Branko Malešević, associate professor at Faculty of Electrical Engineering, Department of Applied Mathematics, University of Belgrade, Serbia, e-mail: malesevic@etf.rs

[3] MSc Bojan Banjac, student of doctoral studies of Software Engineering, Faculty of Electrical Engineering, University of Belgrade, Serbia,
assistant at Faculty of Technical Sciences, Computer Graphics Chair, University of Novi Sad, Serbia, e-mail: bojan.banjac@uns.ac.rs

[4] PhD Ratko Obradović, full professor at Faculty of Technical Sciences, Computer Graphics Chair, University of Novi Sad, Serbia, e-mail: obrad_r@uns.ac.rs

Minkowski metrics contains taxicab metric for value $k = 1$ and Euclidean metric for $k = 2$. The term „taxicab geometry" was first used by K. Menger in the book [9]. Distances between two points $A$ and $B$ are presented in the Figure 1: $d_1$ (short/lond dished lines) in taxicab metric and $d_2$ (continuous line) in Euclidean metric.

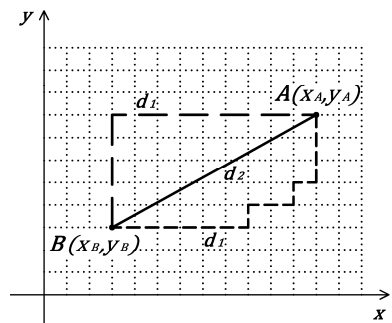

Figure 1. Minkowski distances

Mathematicians Kaya et al. in the paper [6] determined distance from point to line in taxicab geometry with the following statement:

*Lemma 1.* Distance of point $P(x_P, y_P)$ to the line $l: ax + by + c = 0$ in the taxicab plane is:

$$d_1(P, l) = \frac{|ax_P + by_P + c|}{max(|a|,|b|)}. \quad (2)$$

In the same paper, Kaya et al. introduce general equation of taxicab conics with the following statement:

*Theorem 1.* Equation of a taxicab conic, with the foci $(x_1, y_1)$ and $(x_2, y_2)$ or with the focus $(x_1, y_1)$ and directrix $ax + by + c = 0$, has the following form:

$$|x - x_1| + |y - y_1| + \alpha(|x - x_2| + |y - y_2|) + \beta(|ax + by + c|) \mp \alpha\gamma = 0 \quad (3)$$

where $\alpha \in \{-1,0,1\}$, $\beta = \frac{e(\alpha^2 - 1)}{max(|a|,|b|)}$, $\gamma \leq 0$ and $e$ is the eccentricity of related conic.

## 1.1 The Conics in Taxicab Geometry

In this section we give visual presentation of taxicab conics:

*A) Circle* is locus of points with constant distance from one focus. Taxicab circle is square with sides oriented at a 45° angle to the coordinate axes (*Figure 2*). Equation of taxicab circle with focus (centre) $C(x_1, y_1)$ and radius $r$ is given by:

$$|x - x_1| + |y - y_1| = r \qquad (4)$$

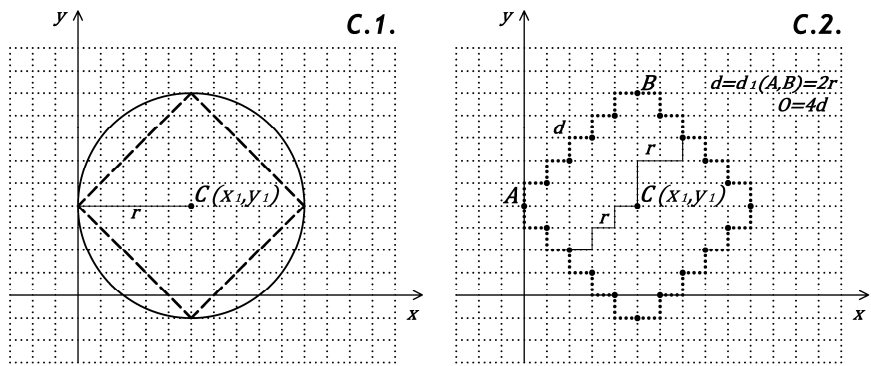

*Figure 2. Circles*

*C.1. Continuous taxicab geometry (dished lines), Euclidean geometry (continuous line);*

*C.2. Discrete taxicab geometry (dots).*

Length of side of square is $r\sqrt{2}$ in Euclidean geometry, while in taxicab geometry this distance is $2r$. In this geometry perimeter of the circle is $8r$, while its area is $4r^2$. Replacement for number $\pi$ in taxicab geometry is number 4.

*B) Ellipse* is locus of points whose sum of distances to two foci is constant. It is polygon in taxicab geometry (*Figure 3. E.1-3*) based on the equation:

$$|x - x_1| + |y - y_1| + |x - x_2| + |y - y_2| + \gamma = 0. \qquad (5)$$

Equation of taxicab ellipse is determined with foci $F_1(x_1, y_1)$, $F_2(x_2, y_2)$ and constant $\gamma \leq 0$. Complete classification of types of polygon that defines taxicab ellipse is given in *Table 1.* from paper [6] based on comparison of parameter $\gamma$ and parameter $\delta = d_1(F_1, F_2)$.

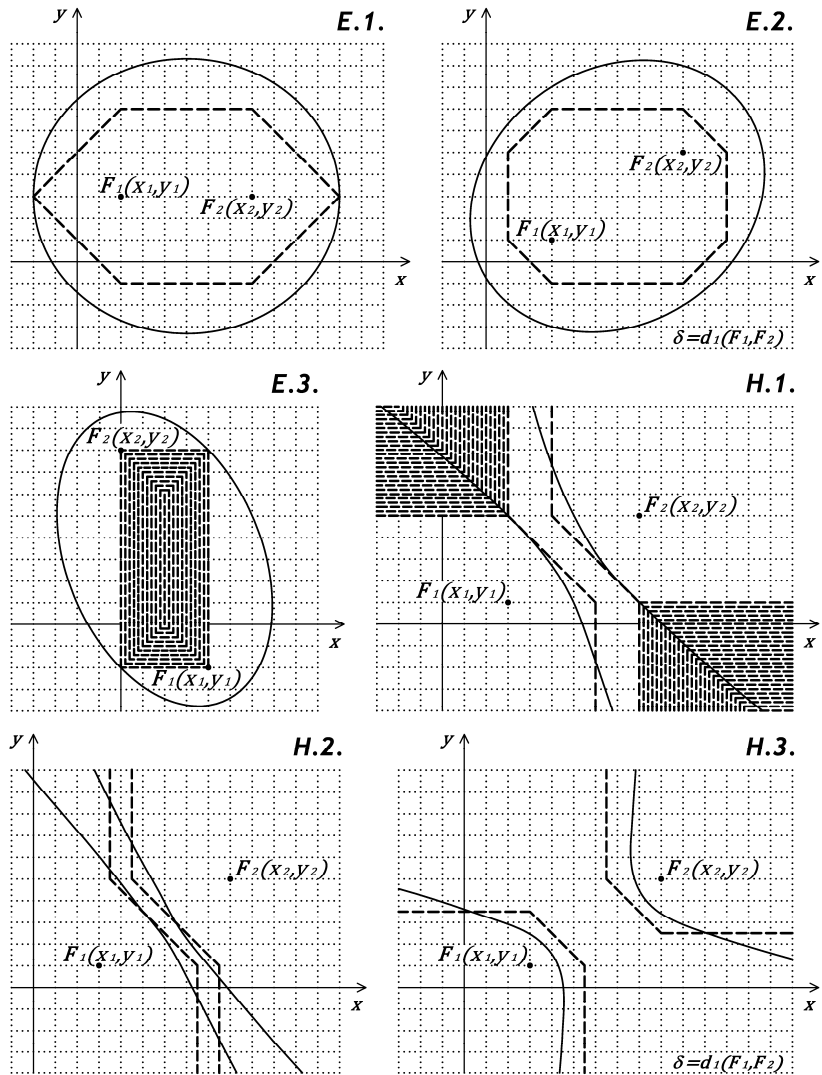

*Figure 3. Ellipses and Hyperbolas*

E.1. If $-\gamma > \delta$ and ($y_1 = y_2$ or $x_1 = x_2$) then ellipse in taxicab geometry become a hexagon;

E.2. If $-\gamma > \delta$ and ($y_1 \neq y_2$ and $x_1 \neq x_2$) then taxicab ellipse is a octagon;

E.3. For $-\gamma = \delta$, ellipse in taxicab geometry become a rectangular region with diagonal $F_1F_2$;

*H.1.* For $-\gamma = \pm(x_1 - x_2 - y_1 + y_2)$, hyperbola in taxicab geometry become two planar regions with tails. Each tail consists of a line segment and vertical/horizontal ray;

*H.2.* If $-\gamma < \delta$ and $-\gamma < |x_1 - x_2 - y_1 + y_2|$ then taxicab hyperbola is a pair of parallel degenerate lines;

*H.3.* If $|x_1 - x_2 - y_1 + y_2| < -\gamma < \delta$ then two-foci taxicab hyperbola become a true taxicab hyperbola.

The perimeter of an ellipse in Euclidean metric can't be exactly calculated. There are several methods which give approximate solution of perimeter of the Euclidean ellipse. G. Almkvist and B. Berndt in the paper [1, pp. 600-601] give a larger number of formulas for approximation of perimeter with error approximations. The area of ellipse in this metric is given by formulae $\pi ab$, where $a$ and $b$ are axis of the ellipse.

Let us emphasize $2\alpha = -\gamma - \delta$. The area of taxicab ellipse *E.1.* (*Figure 3*) is given by $2\alpha \cdot (-\gamma)$, and its perimeter is $2(-\gamma + 2\alpha)$.

The area of taxicab ellipse *E.2.* (*Figure 3*) is given by formula $(|x_1 - x_2| + 2\alpha) \cdot (|y_1 - y_2| + 2\alpha)$. The perimeter of this ellipse is $2(|x_1 - x_2| + 2\alpha) + 2(|y_1 - y_2| + 2\alpha)$.

The area of taxicab ellipse *E.3.* (*Figure 3*) is given by formula $|x_1 - x_2| \cdot |y_1 - y_2|$, and its perimeter is $2\delta$.

*C) Hyperbola* is locus of points whose difference of distance from two foci is constant. It becomes broken line in taxicab geometry (*Figure 3. H.1-3*) based on equation:

$$|x - x_1| + |y - y_1| - (|x - x_2| + |y - y_2|) \mp \gamma = 0. \qquad (6)$$

Equation of taxicab hyperbola is determined by foci $F_1(x_1, y_1)$, $F_2(x_2, y_2)$ and constant $\gamma \leq 0$. In table 2 Kaya et al. give complete classification types of taxicab hyperbolas based on comparison of parameter $\gamma$ and parameter $\delta = d_1(F_1, F_2)$.

*D) Parabola* is locus of points where distance from the foci $F_1(x_1, y_1)$ is equal to $e$-distance from the directrix $d: ax + by + c = 0$. The $e$-distance is standard distance multiplied by $e$-eccentricity of the conic [6]. It becomes broken line in taxicab geometry (*Figure 4. P.1-6*) based on equation:

$$|x - x_1| + |y - y_1| - e(max(|a|, |b|))^{-1}|ax + by + c| = 0. \qquad (7)$$

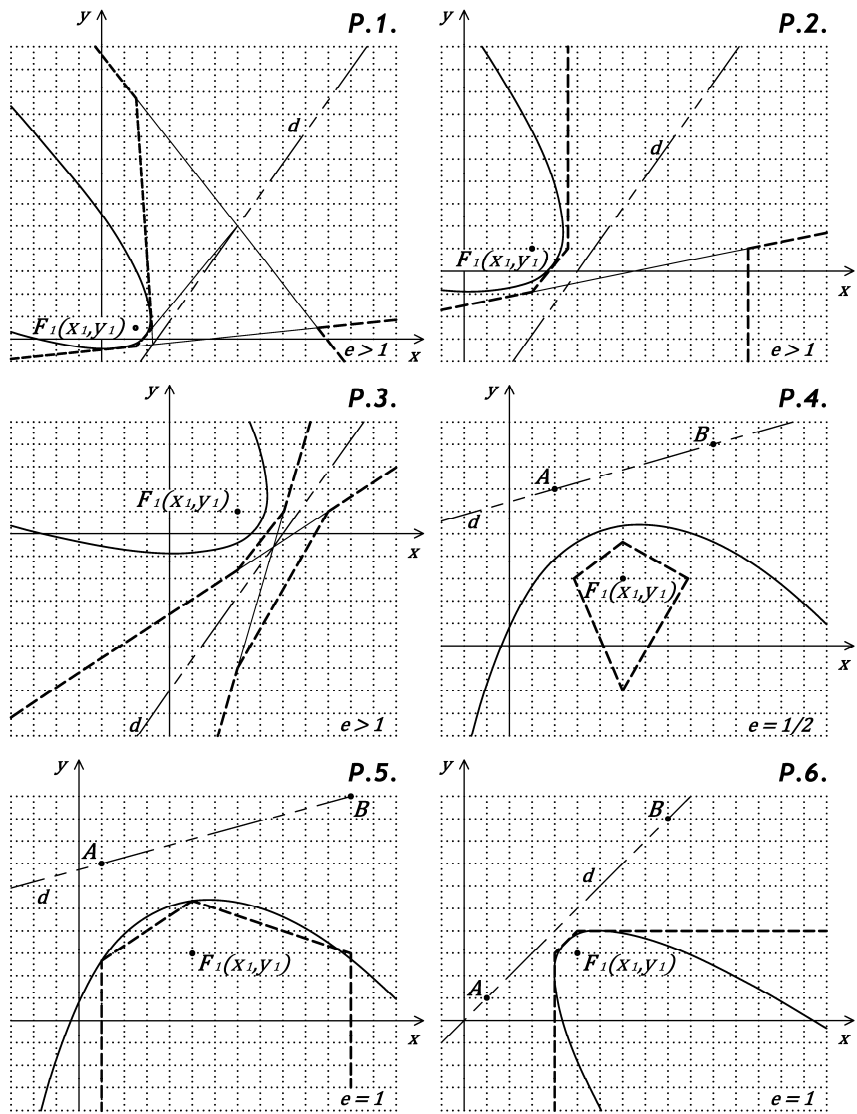

*Figure 4. Parabolas*

*P.1.* If $1 < e < \left|\frac{-a}{b}\right|$ then taxicab parabola consists two branches;

*P.2.* If $1 < e = \left|\frac{-a}{b}\right|$ then taxicab parabola consists two branches;

*P.3.* If $1 < \left|\frac{-a}{b}\right| < e$ then taxicab parabola consists two branches;
6

*P.4.* If $0 < e < 1, \left(\frac{-a}{b} \in \mathbb{R}\right)$ then taxicab parabola become a quadrilateral with a vertical and a horizontal diagonal, both passing through focus;

*P.5.* For $e = 1$ and $\left(\left|\frac{-a}{b}\right| < 1 \text{ or } \left|\frac{-a}{b}\right| > 1\right)$ taxicab parabola become a true taxicab parabola with two line segments and two vertical/horizontal rays;

*P.6.* For $e = 1$ and $\left|\frac{-a}{b}\right| = 1$ taxicab parabola become a true taxicab parabola with a line segment, a vertical ray and a horizontal ray.

## 2. THE TRIFOCAL ELLIPSE IN TAXICAB GEOMETRY

Determining locus of points in plane whose sum of distance to three fixed points is topic which was researched by numerous mathematicians. Historically observed, that researches all have roots from geometric construction of focal curves that appear in papers of R. Descartes [2] as well as J. C. Maxwell [8]. Newer research come from P. V. Sahadevan [11], [12] which define term of egglipse as curve form which sum of distance to tree foci is constant. Sahadevan defines:
a) Egglipse is normal, if foci are collinear and one of foci is not congruent with centre of line segment whose end points are two other foci;
b) Egglipse is abnormal, if foci are collinear and one of foci is in centre of line segment whose end points are two other foci.

By analysing normal trifocal ellipse, Sahadevan had found parametric equations of this curve, and base on them calculated area by use of elliptic integral of first kind. Based on derivate of Sahadevan's parameterization, method for forming of drainage canals with transversal profile in shape of egglipse is patented in India. Similar parameterization of trifocal ellipse like Sahadevan's in special positions of foci were discussed by A. Varga and C. Vincze [13].

In this paper shall be presented trifocal ellipses in taxicab geometry which are defined by following formula:

$$R_1 + R_2 + R_3 = S \qquad (8)$$

where $R_i = d_1(M, F_i)$, $i = 1..3$, are taxicab distances to foci $F_1(x_1, y_1)$, $F_2(x_2, y_2)$, $F_3(x_3, y_3)$ and $S$ given constant. For $S > S_0$ trifocal ellipse is closed curve of egg-like shape (*Figure 5*). For value $S = S_0$ trifocal sum $R_1 + R_2 + R_3$ degenerates to Fermat-Torricelli point [3].

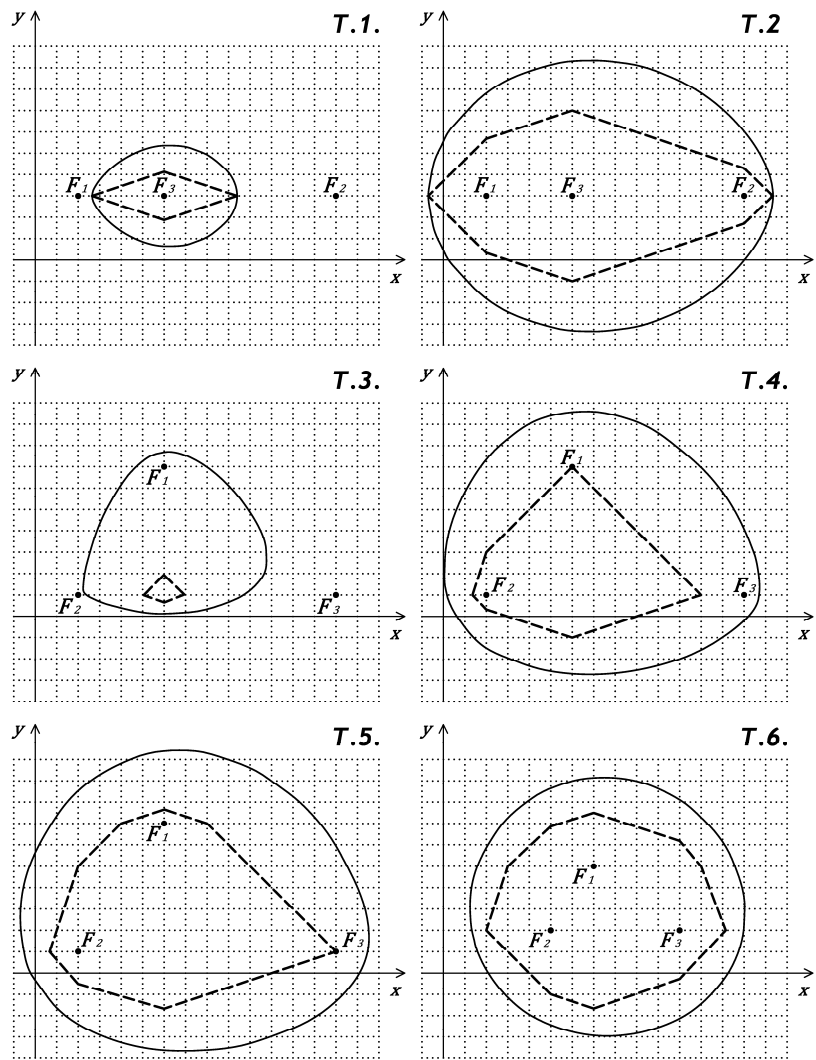

*Figure 5. Trifocal ellipses*

In the paper [13] A. Varga and C. Vincze had considered area $A_E$ and perimeter $P_E$ of trifocal ellipse in Euclidean metric when foci are in special position $F_1(-1,0)$, $F_2(1,0)$ and $F_3(0,0)$. By use of computer-algebra system MAPLE, they got following results:

- for value $S = \frac{5}{2} \Rightarrow A_E = 0.5645$ and $P_E = 2.7123$;
- for value $S = 3 \Rightarrow A_E = 1.7758$ and $P_E = 4.9603$;
- for value $S = 4 \Rightarrow A_E = 4.4032$ and $P_E = 7.5085$.

If we observe taxicab trifocal ellipse with these foci, then area is

$$A_T = \begin{cases} \frac{4}{3}(S-2)^2, & 2 < S < 3 \\ \frac{4}{3}\left(S\left(\frac{S}{3}-1\right)+1\right), & S \geq 3 \end{cases} \quad (9)$$

and perimeter is

$$P_T = \begin{cases} \frac{16}{3}(S-2), & 2 < S < 3 \\ \frac{8}{3}(S-1), & S \geq 3 \end{cases}. \quad (10)$$

For value $S_0 = 2$ trifocal sum $R_1 + R_2 + R_3$ degenerates to Fermat-Torricelli point. For values from the paper [13] ($S > S_0$) we got:

- $S = \frac{5}{2} \Rightarrow A_T = \frac{1}{3} = 0.333\dot{3}$ and $P_T = \frac{8}{3} = 2.666\dot{6}$;
- $S = 3 \Rightarrow A_T = \frac{4}{3} = 1.333\dot{3}$ and $P_T = \frac{16}{3} = 5.333\dot{3}$;
- $S = 4 \Rightarrow A_T = \frac{28}{9} = 3.111\dot{1}$ and $P_T = 8$.

General algorithm for numerical calculation of the area and the perimeter of trifocal ellipses in Euclidean or taxicab metric, implemented in Java application [14], is given by following pseudo-code:

```
AreaAndPerimeter Algorithm
input: startX, endX, step
output: area a, perimeter p
begin
x=startX, a = 0, p = 0
while(x ≤ endX)
  begin
    minY=findMinY(x),  maxY=findMaxY(x)
    if(x = startX) p = p + maxY-minY
      else
       begin
         p = p +  dist (maxY,oldMaxY,step) + dist (minY,oldMinY,step)
         a = a + (maxY-minY)*step  +  t(maxY, oldMaxY, step)
               + t(oldMinY, minY, step)
       end
    oldMinY=minY
    oldMaxY=maxY
    x = x + step
    if ( x > endX) p = p + maxY - minY
  endWhile
end
```

*t algorithm*
input: $Y_1$, $Y_2$, step
output: area a
begin
$$a = (Y_1 - Y_2) * \frac{step}{2}$$
end

*dist algorithm*
input: Y, oldY, step
output: distance d
begin
  if(metrika=TaxiCab) $d = |Y - oldY| + |step|$
  else $d = \sqrt{(Y - oldY)^2 + step^2}$
end

    Algorithm relies on functions *dist* that calculates distance between two points. Functions findMinY and findMaxY seek minimal and maximal values of y for given value of x where function $f(x, y) = R_1 + R_2 + R_3 + S$ gives non positive value. Algorithm needs as starting points minimal and maximal value of x that are part of trifocal ellipse. Its precision depends of numerical precision of functions findMinY and findMaxY, and size of parameter step. For Euclidean metric our experience, based on large number of test cases, is that algorithm makes mistake in calculation of parameter of less than 2% for step value of 1, and lesser for value of area. For taxicab metric error is less than 0.1% .

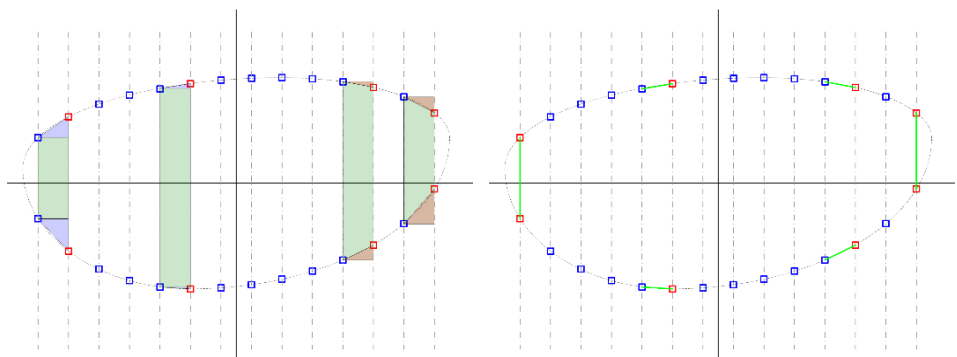

*Figure 6. Visualization of algorithm for area (left) and perimeter (right)*

## 3. CONCLUSION

Specific application of trifocal curves in Euclidean metric is discussed in the paper [10]. Geometry of these curves was applied on solving problem of optimal location in architecture, urbanism and spatial planning. In mentioned paper was conducted concurrent analysis of suggested solution for optimal location Fermat-Torricelli type for existing solution of infrastructural corridors [10] and solution which is obtained using Java applet (available at [14]). Further development of our research shall be based on application of trifocal curves in taxicab metric in architecture, urbanism and spatial planning.

*Acknowledgement.* Research is partially supported by the Ministry of Science and Education of the Republic of Serbia, Grant No. III 44006 and ON 174032